\documentclass[12pt]{article}
 \usepackage[german]{babel}
\usepackage{amssymb}
\usepackage{amsmath}
\usepackage{graphicx}
\usepackage{bbm}
\usepackage[right]{eurosym}
\usepackage[left=2cm,right=2cm,top=2cm,bottom=3cm]{geometry}
\def\vs{\vspace}
\def\noi{\noindent}

\def\IN{\mathbb N}

\def\IR{\mathbb R}
\def\IC{\mathbb C}

\def\IH{\mathbb H}
\def\IE{\mathbb E}

\def\exp{\mathrm{exp}}
\def\sph{\sphericalangle}

\def\Re{\mathrm{Re}}
\def\Im{\mathrm{Im}}
\def\arg{\mathrm{arg}}
\def\ord{\mathrm{ord}}

\begin{document}
\begin{center}
{\bf \Large Asymptotic behaviour of the Riemann mapping function at analytic cusps}
\end{center}

\centerline{Tobias Kaiser and Sabrina Lehner}

\vspace{0.7cm}\noi \footnotesize {{\bf Abstract.} We completely describe the asymptotic behaviour of the Riemann mapping function and its derivatives at an analytic cusp. We achieve the same for the inverse of the mapping function.

\normalsize
\section{Introduction and results}

\subsection{Introduction}
In the literature, there exists a complete asymptotic description of the mapping function at an analytic corner (i.e., the boundary at a given boundary point - let's say the origin - of a given simply connected and proper domain in the complex plane consists of two regular analytic arcs with non-zero opening angle).
The work of Lichtenstein and Warschawski [4, 7, 9] establishes the asymptotic behaviour of the mapping function at an analytic corner and of its derivatives. The mapping function behaves like $z^{1/\alpha}$, where $\pi\alpha$ is the opening angle of the analytic corner (we consider mapping functions onto the upper half plane and assume, applying Carath\'{e}odory's theorem on prime ends, that the origin is mapped to the origin).
Lehman [3] has then obtained the stronger result that the mapping function at an analytic corner (and its inverese) can be developed in a generalized power series. This expansion completely characterizes the asymptotics of the mapping function. We also refer to Pommerenke [5, Chapter 3].

The latter result was used in [1] to embed the Riemann mapping theorem at analytic corners into modern real geometry. Hence, from a geometric point of view it is desirable and necessary to understand the mapping function also at analytic cusps; i.e. if the angle between the two regular analytic arcs vanishes.
But there is a gap in the literature. There are results for arbitrary cusps (see for example Warschawski [8]), but not suitable ones for the analytic setting (and such  are also not just corollaries from the general ones).

In [2], the first author has given the asymptotic behaviour of the mapping function and its inverse in the special case of so-called analytic cups with small perturbation of angles. Only upper bounds for the derivatives of these mappings have been established.
In the present paper, we are able to obtain in the general case of arbitrary analytic cusps a nice description of the asymptotic behaviour of the mapping function, its inverse and the derivatives of these mappings.
The quantities used in the description are geometric invariants of the given analytic cusp.
The precise results will be presented below after a preliminary section, introducing the setting. In the main part of the article, the results are proven. 

 \vs{0.75cm}
\hrule

\vs{0.4cm}
{\footnotesize{\itshape 2010 Mathematics Subject Classification:} 30C20, 30E15}
\newline
{\footnotesize{\itshape Keywords and phrases:} mapping function, analytic cusp, asymptotic behaviour}
\newline
{\footnotesize{\itshape Acknowledgements:} The authors were supported in parts by DFG KA 3297/1-2.}

\subsection{Preliminaries}

By $\IN=\{1,2,3,\ldots\}$ we denote the set of natural numbers, by $\IR$ and $\IC$ we denote the field of real and complex numbers respectively. We let $\IR^*:=\IR\setminus\{0\}$
and $\IC^*:=\IC\setminus\{0\}$ be the set of non-zero real and complex numbers respectively and $\IR_{>0}:=\{x\in\IR\colon x>0\}$ be the set of positive real
numbers. 
Given $z \in \IC$ we denote by $\Re(z)$ its real and by $\Im(z)$ its imaginary part. Furthermore,
$\vert z \vert$ denotes the Euclidean norm of $z \in \IC$ and $\arg(z) \in ]-\pi, \pi]$ the standard argument of $z$. By $\IH := \{z \in \IC \colon \Im(z) > 0\}$ we denote the upper half plane in $\IC$ and by $\IE:= \{z \in \IC \colon \vert z \vert <1\}$ the open unit disc. Given $a\in\IC$ and $r\in\IR_{>0}$, we set $B(a,r) := \{z \in \mathbb{C} \colon\vert z - a \vert < r\}$ and $\overline{B}(a,r):=\{z\in\IC\colon \vert z-a\vert\leq r\}$. A domain is a non-empty open and connected subset of $\IC$.

\vs{0.2cm}
Let $A \subset \mathbb{C}$ be a nonempty set with $0 \in \overline{A}$. Let $f,g: A \rightarrow \mathbb{C}$ be functions such that $g(z) \neq 0$ for all $z$ in a punctured neighbourhood of $0$. We write
\begin{itemize}
\item[(a)] $f \simeq g$ at $0$ on $A$ if $\lim_{z\rightarrow 0} \frac{f(z)}{g(z)} = 1,$
\item[(b)] $f \sim g$ at $0$ on $A$ if $\lim_{z\rightarrow 0} \frac{f(z)}{g(z)} \in \mathbb{C}^*,$
\item[(c)]  $f = o(g)$ at $0$ on $A$ if $\lim_{z\rightarrow 0} \frac{f(z)}{g(z)} = 0,$
\item[(d)] $f = O(g)$ at $0$ on $A$ if there is some $C >0$ such that $\vert f(z)\vert \leq C\vert g(z)\vert$ for all $z$ in a neighbourhood of $0$.
\end{itemize}

\vs{0.2cm}
Given a non-vanishing real or complex power series $h(t) = \sum_{j=0}^\infty a_jt^j$, its {\bf order} is given by
$\text{ord}(h) := \min\{j \in \mathbb{N}_0: a_j \neq 0\}$
and its {\bf leading coefficient} by
$\mathrm{lc}(h):=a_{\mathrm{ord}(h)}.$
Moreover, we set $\text{ord}(0):= \infty$.

\subsection{Results}

Let $\Omega$ be a proper and simply connected domain in $\mathbb{C}$. After applying a translation, we assume that $0\in\partial \Omega$.

A regular analytic curve in the plane is an analytic mapping $\gamma:I\to \IC$ such that $\gamma'(t)\neq 0$ for all $t\in I$ where $I$ is an interval.
A compact trace of a regular analytic curve is called a regular analytic arc.

\vs{0.5cm}
{\bf 1.1 Definition }

\vs{0.1cm}
We say that $\Omega$ has an \textbf{analytic cusp} at $0$ if the boundary of $\Omega$ at $0$ consists of two regular analytic arcs such that the opening angle of $\Omega$ at $0$ vanishes.

\vs{0.5cm}
Assume from now on that $\Omega$ has an analytic cusp at $0$.
We find $\varepsilon>0$ and analytic curves $\gamma,\widetilde{\gamma}:]-\varepsilon,\varepsilon[\to \IC$ with non-vanishing derivatives such that the boundary of $\Omega$ at $0$ is given by the regular analytic arcs $\Gamma:=\gamma([0,\varepsilon/2])$ and $\widetilde{\Gamma}:=\widetilde{\gamma}([0,\varepsilon/2])$.

We perform the following {\bf reductions}:
By shrinking $\varepsilon>0$ if necessary we obtain after complexification that
$\widetilde{\gamma}:B(0,\varepsilon)\to \IC$ is a conformal map onto its image.
Applying the inverse $\widetilde{\gamma}^{-1}$ we can assume that $\widetilde{\gamma}$ is given by the identity and that $\widetilde{\Gamma}$ is a segment contained in the positive real axis.
In this situation, the arc $\Gamma$ is also tangent to the positive real axis.
We write
$$\gamma(t)=|\gamma(t)|\exp\Big(i\arg\big(\gamma(t)\big)\Big)$$
in polar coordinates. For $t$ sufficiently close to $0$ we have
$$\arg(\gamma(t)) =\arctan\left(\frac{\Im\big(\gamma(t)\big)}{\Re\big(\gamma(t)\big)}\right)=:\eta(t).$$
Since $\gamma(0)=0$ and $\gamma'(0)\in\IR_{>0}$ we have $\ord\big(\Re(\gamma(t))\big)=1$ and $\ord\big(\Im(\gamma(t))\big) \geq 2$. Hence, it follows that $|\gamma(t)|$ and $\eta(t)$ are real analytic at $0$ and that $\ord\big(|\gamma(t)|\big) = 1$. Therefore, $|\gamma(t)|$ is locally invertible at $0$ and we can parameterize $\gamma$ close to $0$ by the distance from the origin using the parameter $s := |\gamma|^{-1}(t)$. We obtain $\gamma(s) = s \exp\big(i\mu(s)\big)$ where $\mu(s)$ is real analytic at $0$. Relabeling $\gamma$ and $\widetilde{\gamma}$, we may assume that $\mu(s)$ is positive for small positive $s$.

This transformation is a geometrical invariant of $\Omega$. From the asymptotic behaviour of the Riemann mapping function of the transformed domain one can compute the asymptotic behaviour of the Riemann mapping function of the original domain.

\vs{0.5cm}
In view of the above reductions we assume from now on the following setup for the rest of the paper.

\vs{0.5cm}
{\bf General assumption}

\vs{0.1cm}
{\it There is some $R>0$ such that
$$\Omega\cap \overline{B}(0,R)=\Big\{z\in\IC\colon |z|\leq R, 0<\arg(z)<\sph_\Omega(|z|)\big\}$$
where $\sph_\Omega(t)$ is a real power series that converges on $]-2R,2R[$ and is positive on $]0,R[$.}

\vs{0.5cm}
In this setting, the two boundary arcs of $\Omega$ at $0$ are given by $\Gamma=\gamma([0,R])$ where $\gamma:]-2R,2R[\to \IC, t\mapsto t\exp\big(i\sph_\Omega(t)\big),$ and $\widetilde{\Gamma}=[0,R]$.

\vs{0.5cm}
{\bf 1.2 Definition}

\vs{0.1cm}
The power series $\sph_\Omega(t)$ is called the {\bf angle function} of $\Omega$. We call $N_\Omega=N:=\mathrm{ord}\big(\sph_\Omega(t)\big)$ the \textbf{order of tangency} of $\Omega$ and $a_\Omega=a:=\mathrm{lc}\big(\sph_\Omega(t)\big)$ the \textbf{coefficient of tangency} of $\Omega$.

\vs{0.5cm}
Let $\sum_{j=N_\Omega}^\infty a_{\Omega,j}t^j=\sum_{j=N}^\infty a_jt^j$ be the power series expansion of $\sph_\Omega(t)$. Note that $a=a_N$.

Let $\sph^{-1}_\Omega(t)$ be the multiplicative inverse of $\sph_\Omega(t)$. We write this Laurent series in the form
$t^{-N_\Omega}\sum\limits_{j=0}^\infty b_{\Omega,j} t^j=t^{-N}\sum\limits_{j=0}^\infty b_j t^j$.
For $j\in\{0,\ldots,N-1\}$ we set
$c_{\Omega,j}=c_j:= \pi b_j/(j-N)$.
Moreover, we set $\sigma^\Omega=\sigma:= \pi b_N$.

\vs{0.5cm}
{\bf 1.3 Definition}

\vs{0.1cm}
The tuple $(c_{\Omega,0},\dots,c_{\Omega, N-1},\sigma_\Omega) \in \mathbb{R}^{N+1}$ is called the \textbf{asymptotic tuple of $\Omega$}.

\vs{0.5cm}
The asymptotic tuple of $\Omega$ is a {\bf geometric invariant} of the cusp of the domain $\Omega$ at $0$.

\vs{0.5cm}
Let $\Phi: \Omega \rightarrow \mathbb{H}$ be a conformal map with $\Phi(0) =0$. Then the arc $\widetilde{\Gamma}$ is mapped to the positive real axis and the arc $\Gamma$ is mapped to the negative real axis.

Let $N$ be the order and $a$ be the coefficient of tangency of $\Omega$. Let $(c_0,\ldots,c_{N-1},\sigma)\in\IR^{N+1}$ be the asymptotic tuple of $\Omega$.

\vs{0.5cm}
Theorem A and B completely describe the asymptotic behaviour of the Riemann mapping $\Phi$ and its derivatives.

\vs{0.5cm}
{\bf Theorem A}

\vs{0.1cm}	
{\it We have
$$\Phi(z) \sim z^\sigma \,\exp\left(\frac{c_0}{z^N} + \frac{c_1}{z^{N-1}} + \dots + \frac{c_{N-1}}{z}\right)$$
at $0$ on $\Omega$.}

\vs{0.5cm}
{\bf Theorem B}

\vs{0.1cm}
{\it For $k\in \IN$ we have that
$$\Phi^{(k)}(z) \sim \Phi(z)\,z^{-k(N+1)} \sim z^{\sigma-k(N+1)}\,\exp\left(\frac{c_0}{z^N} + \dots + \frac{c_{N-1}}{z}\right)$$
at $0$ on $\Omega$.}

\vs{0.5cm}
Let $\Psi:\IH\to \Omega$ be a conformal map with $\Psi(0)=0$.
Theorem C\footnote{Recently, Prokhorov [6] has obtained a similar theorem in a different setting of analytic cusps} und D completely describe the asymptotic behaviour of the Riemann mapping $\Psi$ and its derivatives.

\vs{0.5cm}
{\bf Theorem C}

\vs{0.1cm}
{\it We have
$$\Psi(z) \simeq \left(-\frac{\pi}{aN\log(|z|)}\right)^{\frac{1}{N}}$$
	at $0$ on $\mathbb{H}$.}

\vs{0.5cm}
{\bf Theorem D}

\vs{0.1cm}
{\it For $k \in \mathbb{N}$ we have that
$$\Psi^{(k)}(z) \sim \big(\Psi(z)\big)^{N+1}z^{-k} \sim \left(-\frac{1}{\log(z)}\right)^{1+\frac{1}{N}}z^{-k}$$
at $0$ on $\mathbb{H}$.}

\section{Proof of the results}

We use the above notations:
Let $N$ be the order and $a$ be the coefficient of tangency of $\Omega$. Let $(c_0,\ldots,c_{N-1},\sigma)\in\IR^{N+1}$ be the asymptotic tuple of $\Omega$.
Let $\sph_\Omega(t)=\sum_{j=N}^\infty a_jt^j$ and let $\sph_\Omega^{-1}(t)=t^{-N}\sum_{j=0}^\infty b_jt^j$.

Let $\Phi:\Omega\to \IH$ and $\Psi:\IH\to \Omega$ be conformal maps which map $0$ to $0$.

\subsection{Proof of Theorem A}

{\bf 2.1 Lemma}
	
\vs{0.1cm}
{\it Consider the M\"obius transformation
$f:\IH\to \IE_1, z\mapsto 2z/(z+i),$
where $\IE_1:=\{z\in\IC\colon |z-1|<1\}$.
Let $\varphi:=f\circ\Phi$.
Then
$\Phi(z) \sim \varphi(z)$
and
$\arg\big(\Phi(z)\big) = \arg\big(\varphi(z)\big) + \frac{\pi}{2} + o(1)$
at $0$ on $\Omega$.}
	
\vs{0.1cm}
{\bf Proof:}	

\vs{0.1cm}
Let $\sum_{j=1}^\infty d_jz^j$ be the power series expansion of $f$ at $0$. Then $d_1\neq 0$.
From $\varphi(z)=\sum_{j=1}^\infty d_j\big(\Phi(z)\big)^j$ we see that
$\Phi(z)\sim \varphi(z)$.
Moreover, we have
$$\arg\big(\varphi(z)\big) = \arg\left(\frac{2\Phi(z)}{\Phi(z)+i}\right)=\arg\big(\Phi(z)\big)-\arg\big(\Phi(z)+i\big).$$
Since
$\lim_{z \rightarrow 0} \arg\big(\Phi(z)+i\big) = \frac{\pi}{2}$
we get
$\arg\big(\Phi(z)\big) = \arg\big(\varphi(z)\big) + \frac{\pi}{2} + o(1).$
\hfill$\Box$

\vs{0.5cm}
{\bf 2.2 Proposition}

\vs{0.1cm}
{\it The following holds:
\begin{itemize}
\item[(1)]
$$|\Phi(z)|\sim |z|^\sigma \,\exp\left( \frac{c_0}{|z|^N} + \frac{c_1}{|z|^{N-1}} + \dots + \frac{c_{N-1}}{|z|}\right),$$
\item[(2)]
$$\arg\big(\Phi(z)\big)=
\pi\arg(z)|z|^{-N}\left(b_0+b_1|z|+\ldots+b_d|z|^d\right)+o(1).$$
\end{itemize}}

\vs{0.1cm}
{\bf Proof:}

\vs{0.2cm}
{\bf Case 1:} $N > 1$.

\vs{0.1cm}
Since $\sphericalangle_\Omega(t)=a_Nt^N+a_{N+1}t^{N+1}+\ldots$ it follows that $\sphericalangle_\Omega''(t) \sim t^{N-2}, \big(\sphericalangle'_\Omega(t)\big)^2/\sphericalangle_\Omega(t) \sim t^{N-2}$
and $t\big(\sphericalangle'_\Omega(t)\big)^2/\sphericalangle_\Omega(t) \sim t^{N-1}$.
Therefore, the integrals
$$\int_0^\delta \sphericalangle''_\Omega(t)dt,\;\;
\int\limits_0^\delta \frac{(\sphericalangle'_\Omega(t))^2}
{\sphericalangle_\Omega(t)}dt,\;\;
\int\limits_0^\delta \frac{t}{\sphericalangle_\Omega(t)}(\sphericalangle'_\Omega(t))^2 dt$$
converge for small $\delta>0$.
From this we get by  Warschawski [8, Theorem XI(A) \& Theorem XI(B)] and Lemma 2.1 that
$$|\Phi(z)|\sim \exp\left( -\pi \int\limits_t^\delta \frac{dr}{r\sphericalangle_\Omega(r)}\right)$$
and
$$\arg\big(\Phi(z)\big) = \pi \frac{\arg(z)-\frac{1}{2}\sph_\Omega(t)}{\sphericalangle_\Omega(t)}+\frac{\pi}{2} + o(1)=\pi\frac{\arg(z)}{\sph_\Omega(t)}+o(1)$$
where $t:=|z|$.
Setting
$h(t) := -\pi\int_t^\delta dr/\big(r\sphericalangle_\Omega(r)\big)$
we obtain for small positive $\delta$
\begin{eqnarray*}
h(t) &=&  -\pi\int\limits_t^\delta \frac{\sphericalangle^{-1}_\Omega(r)}{r}dr
=  -\pi\int\limits_t^\delta \sum_{j=0}^\infty b_j r^{j-(N+1)}dr\\
&=& -\pi \left( \sum_{j=0}^{N-1}b_j\int\limits_t^\delta r^{j-(N+1)}dr + b_N\int\limits_t^\delta r^{-1}dr + \sum_{j=N+1}^\infty b_j\int\limits_t^\delta   r^{j-(N+1)}dr \right)\\
&=& -\pi \left( \sum\limits_{j=0}^{N-1} b_j \left[\frac{1}{j-N} r^{j-N}\right]_t^\delta + \big[ b_N \log(r)\big]_t^\delta +\sum\limits_{j=N+1}^\infty b_j \left[\frac{1}{j-N}r^{j-N}\right]_t^\delta \right)\\
&=& -\pi \left(D - \sum\limits_{j=0}^{N-1} \frac{b_j}{j-N} t^{j-N} - b_N \log(t) - \sum\limits_{j=N+1}^\infty \frac{b_j}{j-N}t^{j-N} \right)
\end{eqnarray*}
with a constant $D \in \mathbb{R}$ (depending on the choice of $\delta$).
Since
$$\lim\limits_{t \rightarrow 0} \sum\limits_{j=N+1}^\infty \frac{b_j}{j-N}t^{j-N} = 0$$
it follows that
$$h(t) = -\pi D + \pi \sum\limits_{j=0}^{N-1} \frac{b_j}{j-N} t^{j-N} + \pi b_N \log(t) + o(1).$$
Therefore, we obtain  that
\begin{eqnarray*}
\vert \Phi(z)\vert &\sim&\exp\left(h(t)\right) \\
&\sim& \exp\left(-\pi D + \pi \sum\limits_{j=0}^{N-1}\frac{b_j}{j-N}t^{j-N} + \pi b_N \log(t)\right)\\
&\sim&\exp\left(\pi \sum\limits_{j=0}^{N-1}\frac{b_j}{j-N}t^{j-N} + \pi b_N \log(t)\right)\\
&=&z^\sigma\exp\left(\frac{c_0}{t^{N}}+\ldots+\frac{c_{N-1}}{t}\right).
\end{eqnarray*}
This shows the desired asymptotic of the modulus. We deal with the argument. We have that
$$\arg\big(\Phi(z)\big) =
\pi \frac{\arg(z)}{\sphericalangle_\Omega(t)} + o(1)\\
= \pi\arg(z)\sum\limits_{j=0}^\infty b_j t^{j-N} + o(1).$$
Since
$\lim_{t\rightarrow 0} \sum\limits_{j=N+1}^\infty b_j t^{j-N} = 0$
we get
$$\arg\big(\Phi(z)\big) = \pi \arg(z) t^{-N}\big(b_0+\ldots+b_N t^N)+ o(1).$$
	
\newpage
{\bf Case 2:} $d=1$.

\vs{0.1cm}
We apply the transformation $\omega: \mathbb{C}\setminus\IR_{\leq 0} \rightarrow \mathbb{C},~ \omega(z) = \sqrt{z},$ to $\Omega$ and set $\widetilde\Omega := \omega(\Omega)$. Then $\widetilde{\Omega}$ fulfils the general assumption and $\sphericalangle_{\widetilde\Omega}(t) = \frac{1}{2} \sphericalangle_\Omega(t^2)$. Therefore,
$\sphericalangle_{\widetilde\Omega}(s) = \sum_{j=\widetilde{N}}^\infty \widetilde a_j t^j$
where $\widetilde{N} = 2$. Let $\widetilde\Phi: \widetilde\Omega \rightarrow \mathbb{H}, z\mapsto \Phi(z^2)$. Since $\sphericalangle_{\widetilde\Omega}$ fulfils the condition of Case 1, we obtain the following estimates for the modulus and the argument of the mapping function $\Phi=\widetilde\Phi\circ\omega$, setting $s=\sqrt{t}$:
$$\vert \Phi(z)\vert = \vert \widetilde{\Phi}(\omega(z)) \vert
\sim s^{\widetilde{\sigma}}
\exp\left(\frac{\widetilde{c}_0}{s^2}+\frac{\widetilde{c}_1}{s}\right),$$
$$\arg\big(\Phi(z)\big) = \arg\big(\widetilde\varphi(\omega(z))\big)=
\pi\arg\big(\omega(z)\big)s^{-2}
\big(\widetilde{b}_0+\widetilde{b}_1s+\widetilde{b}_2s^2\big)$$
where $(\widetilde{c}_0,\widetilde{c}_1,\widetilde{\sigma})$ is the asymptotic tuple of $\widetilde{\Omega}$ and $\sph^{-1}_{\widetilde{\Omega}}=t^{-2}\sum_{j=0}^\infty \widetilde{b}_jt^j$.
From
$\sphericalangle_{\widetilde\Omega}(t) = \frac{1}{2}\sphericalangle_\Omega(t^2)=\sum_{j=1}^\infty (a_j/2)t^2$
one computes $\widetilde{b}_0 =2b_0,\widetilde{b}_1= 0$ and  $\widetilde{b}_2= 2b_1$
and hence $\widetilde{c}_0=2c_0,\widetilde{c_1}=0$ and $\widetilde{\sigma}=2\sigma$.
Using this and pluging in $s=\sqrt{t}$ we end up with
$$\vert \Phi(z)\vert \sim t^\sigma \exp\left(\frac{c_0}{t} \right)$$
for the modulus.
Since $\arg\big(\omega(z)\big) = \arg\big(\sqrt{z}\big) = \frac{1}{2} \arg(z)$ we obtain
$$\arg\big(\Phi(z)\big)=  \frac{\pi}{2}\arg(z)t^{-1}\left( 2b_0  + 2b_1t \right) + o(1)
=  \pi\arg(z)t^{-1}\left( b_0 + b_1 t\right) + o(1).$$
\hfill$\Box$

\vs{0.5cm}
{\bf 2.3 Lemma}

\vs{0.1cm}
{\it Let $A$ be a nonempty subset of $\mathbb{C}$ such that $0 \in \overline{A}$ and let $f,g: A \to \mathbb{C}$. If $f=g+o(1)$ at $0$ then $\exp(f) \simeq \exp(g)$ at $0$.}

\vs{0.1cm}
{\bf Proof:}

\vs{0.1cm}
We have
\begin{eqnarray*}
f(z) &=& g(z) + o(1)\\
&\Leftrightarrow& f(z)-g(z)= o(1)\\
&\Leftrightarrow& \lim\limits_{z\rightarrow 0}\big(f(z)-g(z)\big)=0 \\
&\Rightarrow& \lim\limits_{z\rightarrow 0} \exp\big(f(z)-g(z)\big) = 1 \\
&\Leftrightarrow& \lim\limits_{z\rightarrow 0}\frac{\exp(f(z))}{\exp(g(z))} =1\\
&\Leftrightarrow& \exp(f(z)) \simeq \exp(g(z)).
\end{eqnarray*}
\hfill$\Box$

\newpage
{\bf Proof of Theorem A:}

\vs{0.1cm}
On $\Omega$ we set
$$H(z):=\sum_{j=0}^{N-1}c_jz^{j-N}+\sigma \log(z).$$

\vs{0.2cm}
{\bf Claim 1:} We find some $D\in \IR$ such that $\mathrm{Re}\big(H(z)\big)=D+\log\big(|\Phi(z)|\big)+o(1)$ at $0$.

\vs{0.1cm}
{\bf Proof of Claim 1:}
Let $t:=\vert z\vert$.
By Proposition 2.2(1) we find some $C\in\IR_{>0}$ such that
$$|\Phi(z)|=C\exp\left(H(|z|+o(1))\right)$$
and hence
\begin{eqnarray*}
\log\big(|\Phi(z)|\big) &=& \log\Big(C\,\exp\big(H(|z|)+o(1)\big)\Big) \\
&=& \log(C) + \log\Big(\exp\big(H(|z|)\big) +o(1))\Big)\\
&=&\log(C)+\log\Big(\exp\big(H(|z|)\big)\big(1+o(1)\big)\Big)\\
&=& \log(C) + H(|z|) + o(1).
\end{eqnarray*}
So it is enough to show that $\Re\big(H(z)\big)=H(|z|)+o(1)$.
We have that
$$\Re\big(H(z)\big) = \sum_{j=0}^{N-1}c_j t^{j-N}\cos\big((j-N)\arg(z)\big) + \sigma\log(t).$$
By the power series expansion of cosine we see that
$$\cos\big((j-N)\arg(z)\big) = 1 + O\big(\arg(z)^2\big).$$
Since $\sphericalangle_\Omega(t) = at^N + O(t^{N+1})$ on $\Omega$ we have $0 \leq \arg(z) \leq at^N+O(t^{N+1})$ and thus
$0 \leq \arg(z)^2 \leq a^2t^{2N}+ O(t^{2N+1})$
for $z \in \Omega$ close to $0$.
Using this we obtain
\begin{eqnarray*}
\lim\limits_{z\rightarrow 0} \big(\Re(H(z)) - H(t)\big)&=&
\lim\limits_{z\rightarrow 0} \left(\sum_{j=0}^{N-1}c_j t^{j-N}\Big(\cos\big((j-N)\arg(z)\big)-1\Big)\right)\\
&=& \lim\limits_{z\rightarrow 0} \sum_{j=0}^{N-1}c_j t^{j-N}\Big(O\big(\arg(z)^2\big)\Big)\\
&=& \lim\limits_{z \rightarrow 0} \sum\limits_{j=0}^{N-1} c_j t^{j-N}\left(a^2t^{2N} + O\big(t^{2N+1}\big)\right)\\
&=&~0
\end{eqnarray*}
and Claim 1 is proven.
\hfill$\Box_{\mathrm{Claim}\,1}$

\vs{0.2cm}
{\bf Claim 2:} We have that $\Im\big(H(z)\big)= \arg\big(\Phi(z)\big)+o(1)$ at $0$.

\vs{0.1cm}
{\bf Proof of Claim 2:}
Let $t:=|z|$. We have
$$\Im\big(H(z)\big) = \sum_{j=0}^{N-1}c_j t^{j-N} \sin\big((j-N)\arg(z)\big) + \sigma\arg(z).$$
By the power series expansion of sine we get
$$\sin\big((j-N)\arg(z)\big) = (j-N)\arg(z) + O\big(\arg(z)^3\big).$$
On $\Omega$ close to the origin we have
$$0 \leq \arg(z)^3 \leq a^3t^{3N}+ O(t^{3N+1}).$$
Applying Proposition 2.2(2) and using the definition of the asymptotic tuple we obtain
\begin{eqnarray*}
&&\lim\limits_{z\rightarrow 0} \Big(\Im\big(H(z)\big) - \arg\big(\Phi(z)\big)\Big)\\
&=& \lim\limits_{z\rightarrow 0}   \left(\sum_{j=0}^{N-1}c_j t^{j-N} \sin\big((j-N)\arg(z)\big) + \sigma \arg(z) - \pi \arg(z) \sum\limits_{j=0}^{N} b_j t^{j-N} + o(1)\right)\\
&=& \lim\limits_{z\rightarrow 0} \left(\sum_{j=0}^{N-1} \frac{\pi b_j}{j-N} t^{j-N} \left((j-N)\arg(z) + O\big(\arg(z)^3\big)\right) +\pi b_N \arg(z)\right. \\
\qquad&& \hspace{1cm}\left.-\pi \arg(z) \sum\limits_{j=0}^{N-1} b_j t^{j-N} - \pi b_N \arg(z) +  o(1)\right) \\
&=& \lim\limits_{z\rightarrow 0} \left(\sum_{j=0}^{N-1} \pi b_j t^{j-N} \left(O\big(\arg(z)^3\big)\right) +  o(1)\right)\\
&=& \lim\limits_{z\rightarrow 0} \left(\sum_{j=0}^{N-1} \pi b_j t^{j-N} (a^3t^{3N}+O\big(t^{3N+1})\big) +  o(1)\right)\\
&=&~0
\end{eqnarray*}
and Claim 2 is proven.
\hfill$\Box_{\mathrm{Claim}\,2}$

\vs{0.2cm}
By Claim 1 and 2 and the definition of the complex logarithm
we obtain that at $0$
\begin{eqnarray*}
H(z)&=&\Re\big(H(z)\big)+i\Im\big(H(z)\big)\\
&=&D+\log\big(|\Phi(z)|\big)+i\arg\big(\Phi(z)\big)+o(1)\\
&=&D+\log\big(\Phi(z)\big)+o(1)
\end{eqnarray*}
for some $D\in \IR$.
Applying Lemma 2.3 we deduce that
$$\Phi(z) \simeq \exp\big(-D + H(z)\big) \sim\exp\big(H(z)\big)$$
and Theorem A is shown.
\hfill$\Box$

\vs{0.5cm}
Recall the definition of small perturbation of angles in [2].

\vs{0.5cm}
{\bf 2.4 Definition} (cf. [2, Definition 3])

\vs{0.1cm}
We say that $\Omega$ has {\bf small perturbation of angles} if $\sph_\Omega(t)=at^N+o(t^{2N})$.

\vs{0.5cm}
The main result of [2] (see [2, Theorem 6 \& Lemma 5]) is included in Theorem A:

\vs{0.5cm}
{\bf 2.5 Corollary}

\vs{0.1cm}
{\it If $\Omega$ has small perturbation of angles then
$$\Phi(z) \sim \exp\left(-\frac{\pi}{a N z^N}\right)$$
at $0$ on $\Omega$.}

\vs{0.1cm}
{\bf Proof:}

\vs{0.1cm}
The coefficients of the multiplicative inverse $\sph^{-1}_\Omega(t)=t^{-N}\sum_{j=0}^\infty b_jt^j$ of $\sph_\Omega(t)=\sum_{j=N}^\infty a_jt^j$ are given by
$$b_j =
\left\{\begin{array}{lll}
\frac{1}{a_N},&&j=0,\\
&\mbox{if}&\\
-\frac{1}{a_N}\sum\limits_{k+l=j}b_k a_{N+l}, &&j \neq 0.
\end{array}\right.$$
Since $a_{N+1} = \ldots = a_{2N} = 0$ we obtain
$b_j=0$ for $0< j \leq N$.
By the definition of the asymptotic tuple we get
$c_0 = -\frac{\pi}{aN}, c_j = 0$ for $j \in \{1,\ldots, N-1\}$, and $\sigma = 0$.
Hence Theorem A implies the assertion.
\hfill$\Box$

\vs{0.5cm}
The following example (cf. [2, Remark 8(i)]) does not have small perturbation of angles but can be handled now:

\vs{0.5cm}
{\bf 2.6 Example}

\vs{0.1cm}
{\it Let
$$\Omega := \left\{z \in \mathbb{C} \colon 0 < \vert z\vert < \frac{1}{2},~ 0 < \arg(z) < \vert z\vert-\vert z\vert^2\right\}.$$
Then
$$\Phi(z)\sim z^\pi\,\exp\left(-\frac{\pi}{z}\right).$$}

\vs{0.1cm}
{\bf Proof:}

\vs{0.1cm}
The angle function of $\Omega$ is given by $\sph_\Omega(t)=t-t^2$. So $N_\Omega=1$. By Theorem A we get
$\Phi(z) \sim z^\sigma\exp\left(c_0/z\right)$
where $c_0=-\pi b_0$ and $\sigma= \pi b_1$ and $t^{-1}(b_0+b_1t+\ldots)$ is the multiplicative inverse of $t-t^2$. One computes $b_0=b_1=1$, and we are done.
\hfill$\Box$

\subsection{Proof of Theorem B}

{\bf 2.7 Proposition}

\vs{0.1cm}
{\it Let
$$F: \mathbb{C}\setminus\IR_{\leq 0} \rightarrow \mathbb{C}, z\mapsto z^\sigma\exp\left(\frac{c_0}{z^N} + \dots + \frac{c_{N-1}}{z}\right).$$
Given $k\in\IN$, we have that
$$F^{(k)}(z) \sim F(z) z^{-k(N+1)}$$
at $0$.}

\newpage
{\bf Proof:}

\vs{0.1cm}
Let
$$H(z) := \sum_{j=0}^{N-1}c_j z^{j-N} + \sigma \log(z).$$
By the formula of Fa\`{a} di Bruno for the derivatives of the composition of two functions and the fact that $\exp^{(k)}(z) = \exp(z)$ for $k \in \mathbb{N}$  we see that
\begin{eqnarray*}
F^{(k)}(z) &=& \frac{d^k}{dz^k}(\exp(H(z))) \\
&=& \sum\limits_{(j_1,\ldots,j_k)\in T_k} \frac{k!}{j_1! \cdot \ldots \cdot j_k!} \exp^{(j_1+ \ldots +j_k)}(H(z)) \prod_{l=1}^k\left(\frac{1}{l!} \frac{d^l}{dz^l} H(z)\right)^{j_l} \\
&=& \sum\limits_{(j_1, \ldots ,j_k)\in T_k} \frac{k!}{j_1!\cdot \ldots \cdot j_k!} \exp(H(z)) \prod_{l=1}^k\left(\frac{1}{l!} \frac{d^l}{dz^l}H(z)\right)^{j_l}
\end{eqnarray*}
where $T_k$ is the set of all k-tupels $(j_1,\ldots,j_k)\in \IN_0^k$ such that $1j_1+2j_2+\ldots+kj_k=k$. Since
$$\frac{1}{l!} \frac{d^l}{dz^l} H(z) = \frac{1}{l!} \frac{d^l}{dz^{l}} \left(\sum\limits_{j=0}^{N-1} c_jz^{j-N} + \sigma\log(z)\right)\sim z^{-(N+l)}$$
for $l \in \mathbb{N}_0$ (note that $c_0=1/a\neq 0$) and since $1j_1+\ldots+kj_k = k$ we obtain
\begin{eqnarray*}
\prod_{l=1}^k\left(\frac{1}{l!} \frac{d^l}{dz^l} H(z)\right)^{j_l} &\sim& \prod_{l=1}^k z^{-(N+l)j_l} \nonumber\\
&=& z^{-(N+1)j_1}\cdot \ldots \cdot z^{-(N+k)j_k} \nonumber \\
&=& z^{-(j_1+\ldots +j_k)N - (1j_1+\ldots +kj_k)} \nonumber \\
&=& z^{-(j_1+\ldots +j_k)N-k}.
\end{eqnarray*}
Given $(j_1,\ldots,j_k)\in T_k$, we have that $j_1+\ldots +j_k \in \{1,\ldots ,k\}$ and that $j_1+ \ldots + j_k = k$ if and only if $(j_1,j_2,\ldots,j_k)=(k,0,\ldots,0)$.
Thus, it follows that
$$F^{(k)}(z)
\sim \exp(H(z))z^{-k(N+1)}= F(z)z^{-k(N+1)}.$$
\hfill$\Box$

\vs{0.5cm}
{\bf Proof of Theorem B:}

\vs{0.1cm}
Let $r,s\in\IR$ with $0<s<r$ such that $\gamma(t)=t\exp\big(i\sph_\Omega(t)\big)$ is injective on $B(0,r)$ and that $B(0,s)\subset \gamma(B(0,r))$. Let
$$\Omega':=\left\{z\in \IC\mid \overline{z}\in\Omega\right\}\mbox{ and } \Omega'':=\left\{z\in\IC\mid z\in B(0,s)\mbox{ and }\gamma(\overline{\gamma^{-1}(z)})\in\Omega\right\}$$
and set $\widehat{\Omega}:=\Big(\Omega \cup \Gamma \cup \widetilde{\Gamma} \cup \Omega' \cup \Omega''\Big)\cap B(0,s).$
By reflection at analytic arcs, we have that $\Phi$ has a holomorphic extension $\widehat{\Phi}$ to $\widehat{\Omega}$ given by
$$ \widehat{\Phi}(z) =
\left\{\begin{array}{lll}
\Phi(z), && z \in \Omega \cup \Gamma \cup \widetilde{\Gamma},\\
\overline{\Phi(\overline{z})}, &\mbox{if}& z \in \Omega',\\
\overline{\Phi\left(\gamma\left(\overline{\gamma^{-1}(z)}\right)\right)}, &&
z \in \Omega''.\\
\end{array}\right.$$
The function $F(z)=\exp\big(\sum_{j=0}^{N-1}c_jz^{j-N}+\sigma\log(z)\big)$ is holomorphic on $\mathbb{C}\setminus\IR_{\leq 0}$ and hence on $\widehat{\Omega}$ (after shrinking $r$ and $s$ if necessary).

\vs{0.2cm}
{\bf Claim 1:} $\widehat{\Phi}(z) \sim F(z)$ at $0$ on $\widehat\Omega$.

\vs{0.1cm}
{\bf Proof of Claim 1:}
By Theorem A we have
$\Phi(z) \sim F(z)$
at $0$ on $\Omega\cup\Gamma\cup\widetilde{\Gamma}$. Since
$\overline{F\left(\overline{z}\right)} = F(z)$
for $z\in \IC\setminus\IR_{\leq 0}$
we obtain that
$\widehat{\Phi}(z) = \overline{\Phi(\overline{z})} \sim F(z)$
at $0$ on $\Omega'$. It remains to show that $\widehat{\Phi}(z) \sim F(z)$ at $0$ on $\Omega''$.
From $\gamma(0)=\gamma^{-1}(0)=0$ we get that
$$\overline{\Phi\left(\gamma\left(\overline{\gamma^{-1}(z)}\right)\right)} \sim \overline{F\left(\gamma\left(\overline{\gamma^{-1}(z)}\right)\right)}.$$
at $0$ on $\Omega''$. Therefore, we have to show that
$$\overline{F\left(\gamma\left(\overline{\gamma^{-1}(z)}\right)\right)} \sim F(z)$$
at $0$ on $\Omega''$. We have
$F(\gamma(z)) = \exp\big(H(\gamma(z)\big)$
where
$$H(z)=\sum\limits_{j=0}^{N-1} c_j z^{j-N} + \sigma \log(z).$$
Since $\gamma(z)= z\exp\left(i\sphericalangle_\Omega(z)\right)$ we see by the power series expansion of the exponential function that
$\big(\gamma(z)\big)^{j-N} = z^{j-N}+ o(1)$
for $j\in\{0,\ldots,N-1\}$.
Applying Lemma 2.3, we obtain
\begin{eqnarray*}
F\big(\gamma(z)\big) &\sim& \exp\left(\sum\limits_{j=0}^{N-1} c_jz^{j-N} + \sigma \log\Big(z\exp\big(i\sphericalangle_\Omega(z)\big)\Big)\right)\\
&\sim&\exp\left(\sum\limits_{j=0}^{N-1} c_jz^{j-N} + \sigma \log(z)+\sigma i\sphericalangle_\Omega(z)\right)\\
&\sim&\exp\big(H(z)\big).
\end{eqnarray*}
Hence $F\big(\gamma(z)\big)\sim F(z)$ at $0$ on $\Omega''$ and thus
$$\overline{F\left(\gamma\left(\overline{\gamma^{-1}(z)}\right)\right)} \sim \overline{F\left(\overline{\gamma^{-1}(z)}\right)}.$$
We also have that $\gamma^{-1}(z)=z+O(z^{N+1})$. So similarly to above we get that
$$F\left(\overline{\gamma^{-1}(z)}\right)\sim F\left(\overline{z}\right).$$
Thus, we see that
$$\overline{F\left(\gamma\left(\overline{\gamma^{-1}(z)}\right)\right)}~ {\sim}~ \overline{F\left(\overline{\gamma^{-1}(z)}\right)}~ {\sim} \overline{F\left(\overline{z}\right)}=F(z)$$
and Claim 1 is proven.
\hfill$\Box_{\mathrm{Claim}\,1}$

\newpage
{\bf Claim 2:} There is some $\rho > 0$ such that $B(z,2\rho\vert z \vert^{N+1}) \subset \widehat\Omega$ for all sufficiently small $z \in \Omega \cup \Gamma \cup \widetilde{\Gamma}$.

\vs{0.1cm}
{\bf Proof of Claim 2:}
Since $\sphericalangle_\Omega(t) \sim t^N$ we see that for $t>0$
$$\text{dist}\big(t, \Gamma)\sim \Im\big(\gamma(t)\big)=t\sin\big(\sph_\Omega(t)\big) \sim t^{N+1}$$
at $0$. Reflecting at the positive real axis, we find some $\rho_1 >0$ such that $B(t, 2\rho_1 t^{N+1}) \subset \widehat\Omega$ for all sufficiently small $t>0$.
Similarly, we have $\text{dist}\big(z,\widetilde\Gamma\big)\sim |z|^{N+1}$ at $0$ on $\Gamma$.
Reflecting at the analytic arc $\Gamma$ (note that $\gamma(z)\sim z$) we find some
$\rho_2>0$ such that $B(z,2\rho_2|z|^{N+1})\subset \widehat{\Omega}$ for all $z\in\Gamma$ sufficiently small. Take $\rho:=\min\{\rho_1,\rho_2\}$.
\hfill$\Box_{\mathrm{Claim}\,2}$

\vs{0.2cm}
For $k\in\IN_0$ let $\rho_k:=\rho/2^k$ and
$$\widehat{\Omega}_k:=\big\{z\in \widehat{\Omega}\colon \mathrm{dist}(z,\Omega)< \rho_k|z|^{N+1}\big\}.$$

\vs{0.2cm}
{\bf Claim 3:} $\widehat{\Phi}^{(k)}(z)\sim F^{(k)}(z)$ at $0$ on $\widehat{\Omega}_k$
for all $k\in \IN_0$.

\vs{0.1cm}
{\bf Proof of Claim 3:} We show it by induction on $k$.

\vs{0.2cm}
$\mathbf{k=0:}$ 
We obtain the base case by Claim 1.

\vs{0.2cm}
$\mathbf{k\to k+1:}$
By the inductive hypothesis there exists a holomorphic function
$h:\widehat{\Omega}_k \rightarrow \mathbb{C}$ such that $h(z) = o(1)$ at $0$ on
$\widehat{\Omega}_k$ and some constant $\alpha\in \mathbb{C}^*$ such that
$\widehat{\Phi}^{(k)}(z) = \alpha F^{(k)}(z) +F^{(k)}(z)h(z)$
on $\widetilde{\Omega}_k$.
Hence
$$\widehat{\Phi}^{(k+1)}(z) =
\alpha F^{(k+1)}(z) + F^{(k+1)}(z)h(z) + F^{(k)}(z)h'(z).$$
By Proposition 2.7 we have that $F^{(k+1)}(z)\sim F^{(k)}(z)/z^{N+1}$.
We show that $h'(z)=o(1/z^{N+1})$ on $\widehat{\Omega}_{k+1}$ and are done.

Let $\eta:=\rho_{k+2}$. By Claim 2 we have that $\overline{B}(z,\eta|z|^{N+1})\subset \widehat{\Omega}_k$ for all $z\in\widehat{\Omega}_{k+1}$.
We obtain by the Cauchy estimate that
\begin{eqnarray*}
\vert h'(z) \vert &\leq& \frac{1}{\eta}  \frac{1}{\vert z \vert^{N+1}}\max_{\vert w-z\vert = \eta \vert z \vert^{N+1}} \vert h(w) \vert\\
&=&O\left(\frac{1}{\vert z \vert^{N+1}}\max_{\vert w-z\vert = \eta \vert z \vert^{N+1}} \vert h(w) \vert \right)\\
&=&o\left(\frac{1}{z^{N+1}}\right)
\end{eqnarray*}
at $0$ on $\widehat{\Omega}_{k+1}$.
\hfill$\Box_{\mathrm{Claim\,3}}$

\vs{0.2cm}
Claim 3, Theorem A and Proposition 2.7 finish the proof of Theorem B.
\hfill$\Box$

\vs{0.5cm}
{\bf 2.8 Example}

\vs{0.1cm}
{\it Let $\Omega$ be as in Example 2.6. Then
$$\Phi^{(k)}(z) \sim z^{\pi-2k}\exp\left(-\frac{\pi}{z}\right)$$
at $0$ on $\Omega$ for every $k\in\IN$.}

\subsection{Proof of Theorem C}

{\bf Proof of Theorem C:}

\vs{0.1cm}
Let $\Phi := \Psi^{-1}$. By Theorem A we have
$$\Phi(w) \sim \exp\left( \sum_{j=0}^{N-1}c_jw^{j-N} + \sigma \log(w)\right)$$
at $0$ on $\Omega$ where $\left(c_0, \dots ,c_{N-1},\sigma\right) \in \mathbb{R}^{N+1}$ is the asymptotic tuple of $\Omega$.
We set
$H(w) := \sum_{j=0}^{N-1}c_jw^{j-N} + \sigma \log(w)$
and obtain
$z = \Phi\big(\Psi(z)\big) \sim \exp\big(H(\Psi(z))\big).$
Applying $\log: \mathbb{H} \rightarrow \mathbb{C}$ we get
$\log(z) \simeq H(\Psi(z)).$
Since $\lim\limits_{z\rightarrow 0} \Psi(z) = 0$ we see that
\begin{eqnarray*}
H\big(\Psi(z)\big) &=& c_0\Psi(z)^{-N} + c_1\Psi(z)^{1-N} + \ldots + c_{N-1}\Psi(z)^{-1} + \sigma\log(\Psi(z))\\
&=& \Psi(z)^{-N}(c_0 + c_1\Psi(z)+\ldots+c_{N-1}\Psi(z)^{N-1} + \sigma\log(\Psi(z))\Psi(z)^N)\\
& \simeq& c_0\Psi(z)^{-N}.
\end{eqnarray*}
Hence,
$\log(z) \simeq c_0\Psi(z)^{-N}.$
We are done since $\log(z) \simeq \log|z|$ on $\mathbb{H}$ and $c_0 = -\frac{\pi}{aN}$ where $N$ is the order and $a$ the coefficient of tangency of $\Omega$.
\hfill$\Box$

\vs{0.5cm}
{\bf 2.9 Example}

\vs{0.1cm}
{\it Let $\Omega$ be as in Example 2.6. We have that
$$\Psi(z) \simeq -\frac{\pi}{\log(z)}$$
at $0$ on $\mathbb{H}$.}

\subsection{Proof of Theorem D}

{\bf 2.10 Proposition}

\vs{0.1cm}
{\it Let
$$G: \mathbb{C}\setminus i\IR_{\leq 0} \rightarrow \mathbb{C}, z\mapsto \left(-\frac{1}{\log(z)}\right)^{\frac{1}{N}}.$$
Given $k\in\IN$, we have that
$$G^{(k)}(z) \sim \left(-\frac{1}{\log(z)}\right)^{\frac{1}{N}+1}\frac{1}{z^k}$$
at $0$.}

\vs{0.1cm}
{\bf Proof:}

\vs{0.1cm}
It follows inductively by the product rule that, given $k\in\IN$, there are $d_{k,1},\ldots,d_{k,k}\in\IR$ with $d_{k,1}\neq 0$
such that
$$G^{(k)}(z) = \left( \sum\limits_{j=1}^{k} d_{k,j} \left(-\frac{1}{\log(z)}\right)^{\frac{1}{N}+j}\right) \frac{1}{z^{k}}.$$
This gives the claim.
\hfill$\Box$

\vs{0.5cm}
{\bf Proof of Theorem D:}

\vs{0.1cm}
By the Schwarz reflection principle, $\Psi(z)$ has a holomorphic extension $\widehat{\Psi}$ to $V:=U \cap \mathbb{C} \setminus i\mathbb{R}_{\leq 0}$ where $U$ is a sufficiently small open neighbourhood of $0$ given by
$$\widehat{\Psi}(z) =
\left\{\begin{array}{lll}
\Psi(z),& & z \in U\cap \overline{\mathbb{H}},\\
\overline{\Psi(\overline{z})}, &\mbox{if}& z \in U \cap -\mathbb{H}\mbox{ and }\Re(z)>0,\\
\gamma\left(\overline{\gamma^{-1}\left(\Psi(\overline{z})\right)}\right), && z \in U \cap -\mathbb{H}\mbox{ and }\Re(z)<0.\\
\end{array}\right.$$
Let $G$ be the function of Proposition 2.10.

\vs{0.2cm}
{\bf Claim 1:} $\widehat{\Psi}(z)\sim G(z)$ at $0$ on $V$.

\vs{0.1cm}
{\bf Proof of Claim 1:}
This follows the same lines as the proof of Claim 1 in the proof of Theorem B.
\hfill$\Box_{\mathrm{Claim\,1}}$

\vs{0.2cm}
Let $\widehat{\Phi}:\widehat{\Omega}\to \IC$ be the holomorphic extension function introduced in the proof of Theorem B.
Then $\widehat{\Phi}(z)\sim F(z)$ by Claim 1 in the proof of Theorem B. Hence we find some $C\in \IR_{>0}$ such that $\widehat{\Phi}(\widehat{\Omega}_1)$ contains
$V_1$ where
$$V_k:=\big\{z\in V\colon \mathrm{dist}(z, i\IR_{\leq 0})\geq C^k|z|\big\}$$
for $k\in\IN_0$.

\vs{0.2cm}
{\bf Claim 2:} $\widehat{\Psi}^{(k)}(z)\sim G^{(k)}(z)$ at $0$ on $V_k$ for all $k\in\IN$.

\vs{0.1cm}
{\bf Proof of Claim 2:}
We show it by induction on $k$.

\vs{0.2cm}
$\mathbf{k=1:}$ 
By Claim 1, 
$\widehat{\Psi}(z) \sim G(z)$.
The derivative of $\widehat{\Psi}=\widehat{\Phi}^{-1}$ on $V_1$ can be computed by the inverse function rule as
$\widehat{\Psi}'(z) = 1/\Phi'(\Psi(z))$
on $V_1$. By Claim 3 in the proof of Theorem B we have
$\Phi'(z) \sim \Phi(z) z^{-(N+1)}$ at $0$ on $\widehat{\Omega}_1$.
Hence we obtain by Proposition 2.10 that
$$\widehat{\Psi}'(z)\sim \frac{1}{\Phi\big(\Psi(z)\big)\Psi(z)^{-(N+1)}} = \frac{\Psi(z)^{N+1}}{z}\sim \frac{G(z)^{N+1}}{z}\sim G'(z)$$ 
at $0$ on $V_1$.

\vs{0.2cm}
$\mathbf{k\to k+1:}$
One can copy the inductive step of Claim 3 in the proof of Theorem B, taking into account that $G^{(k+1)}(z)\sim G^{(k)}(z)/z$ at $0$ on $V$ by Proposition 2.10.
\hfill$\Box_{\mathrm{Claim\,2}}$

\vs{0.2cm}
Claim 2, Theorem C and Proposition 2.10 finish the proof of Theorem D.
\hfill$\Box$

\vs{0.5cm}
{\bf 2.11 Example}

\vs{0.1cm}
{\it Let $\Omega$ be as in Example 2.6. Then
$$\Psi^{(k)}(z) \sim \left(-\frac{1}{\log(z)}\right)^2 z^{-k}$$
at $0$ on $\mathbb{H}$ for every $k\in\IN$.}

\newpage
\noi \footnotesize{\centerline{\bf References}

\begin{itemize}
\item[(1)] T. Kaiser: The Riemann Mapping Theorem for Semianalytic Domains
and O-minimality. {\it Proceedings of the London Mathematical Society} (3)
{\bf 98}, no. 2 (2009), 427-444.
\item[(2)] T. Kaiser: Asymptotic Behaviour of the Mapping Function at an Analytic
Cusp with Small Perturbation of Angles. {\it Computational Methods
and Function Theory} {\bf 10} (2010), 35-47.
\item[(3)] R. S. Lehman: Development of the Mapping Function at an Analytic
Corner. {\it Pacific Journal of Mathematics} {\bf 7}, no. 3 (1957), 1437-1449.
\item[(4)] L. Lichtenstein: \"Uber die Konforme Abbildung Ebener Analytischer Gebiete
mit Ecken. {\it Journal f\"ur die reine und angewandte Mathematik} {\bf 140}
(1911), 100-119.
\item[(5)] Ch. Pommerenke: Boundary Behaviour of Conformal Maps, Springer, 1991.
\item[(6)] D. Prokhorov: Conformal Mapping Asymptotics at a cusp. arXiv:1511.00514 (2015), 9 p.
 \item[(7)] S. Warschwaski: \"Uber das Randverhalten der Ableitung der Abbildungsfunktion bei konformer Abbildung. {\it Mathematische Zeitschrift} {\bf 35}
(1932), 321-456.
\item[(8)] S.Warschawski: On Conformal Mapping of Infinite Strips. {\it Transactions
of the American Mathematical Society} {\bf 51}, no. 2 (1942), 280-335.
\item[(9)] S. Warschawski: On a Theorem of L. Lichtenstein. {\it Pacific Journal of
Mathematics} {\bf 5}, no. 5 (1955), 835-839.
\end{itemize}

\vs{0.5cm}
Tobias Kaiser\\
University of Passau\\
Faculty of Computer Science and Mathematics\\
tobias.kaiser@uni-passau.de\\
D-94030 Germany 

\vs{0.2cm}
Sabrina Lehner\\
University of Passau\\
Faculty of Computer Science and Mathematics\\
sabrina.lehner@uni-passau.de\\
D-94030 Germany}
\end{document}